\documentclass[11pt, a4paper,twoside]{article}
\usepackage{amssymb,amsmath,amsthm,latexsym,bbm,mathrsfs,nicefrac,color,xcolor}
\usepackage{a4wide,microtype} % Verschšnerung des Erscheinungsbildes
\usepackage{hyperref} % Links im pdf-Dokument 
\usepackage{natbib} % bibliography ?
\usepackage{graphicx} % Einbinden von Bildern
\bibpunct{(}{)}{;}{a}{,}{,} % format

\newtheorem{theo}{Theorem}[section]
\newtheorem{prop}[theo]{Proposition}
\newtheorem{defi}[theo]{Definition}
\newtheorem{lemm}[theo]{Lemma}

\newtheorem{coro}[theo]{Corollary}
\newtheorem{cons}[theo]{Construction}

\title{The Shark Fin Function --  Asymptotic Behavior of the Filtered Derivative for Point Processes in Case of Change Points}
\author{
{\sc Michael Messer, Gaby Schneider}\\[3ex]
	{Institute of Mathematics}\\{Johann Wolfgang Goethe University}\\[0.5ex]{Robert-Mayer-Str. 10}\\{60325 Frankfurt (Main), Germany}\\[1.5ex]{corresponding author: Michael Messer, messer@math.uni-frankfurt.de}\\
}
\date{}

\begin{document}
	
%\onehalfspacing
\maketitle

\begin{abstract}
A multiple filter test (MFT) for the analysis and detection of rate change points in point processes on the line has been proposed recently. The underlying statistical test investigates the null hypothesis of constant rate. For that purpose, multiple filtered derivative processes are observed simultaneously. Under the null hypothesis, each process $G$ asymptotically takes the form 
\begin{align*}
G \sim  L,
\end{align*}
while $L$ is a  zero-mean Gaussian process with unit variance. This result is used to derive a rejection threshold for statistical hypothesis testing.\\ 
The purpose of this paper is to describe the behavior of $G$ under the alternative hypothesis of rate changes and potential simultaneous variance changes. We derive the approximation
\begin{align*}
G \sim  \Delta \cdot\left(\Lambda + L\right),
\end{align*}
with deterministic functions $\Delta$ and $\Lambda$.
The function $\Lambda$ accounts for the systematic deviation of $G$ in the neighborhood of a change point. When only the rate changes, $\Lambda$ is hat shaped. 
When also the variance changes, $\Lambda$ takes the form of a shark's fin.
% takes the form of a shark's fin whose magnitude is proportional to  a scaled rate difference and grows with the bandwidth of $G$.
In addition, the parameter estimates required in practical application are not consistent in the neighborhood of a change point. Therefore, we derive the factor $\Delta$ termed here the distortion function. It accounts for the lack in consistency and describes the local parameter estimating process relative to the true scaling of the filtered derivative process. 

%Since asymptotics will be given by an increase in the bandwidth, the main impact on the structure of $G$ is induced by the shark fin function $\Lambda$ while $\Delta$ typically plays a minor role in practice. 
 
\end{abstract}\quad

\section{Introduction}
The statistical theory of change point detection aims at the detection of structural breaks (so called change points) in time series. For an overview of the topic see the textbooks of \citet{Brodsky1993, Basseville1993,Csorgo1998} or the review article of \citet{Aue2013}.
 We focus here on renewal processes on the positive line (e.g., \citet{Gut2002,Gut2009,Timmermann2014}). In applications such as neuronal spike trains, structural breaks can occur on different time scales. Interesting multi scale methods have been proposed by \citet{Frick2014,Fryzlewicz2014,Matteson2014}. Particularly for the scenario of  point processes a multiple filter test (MFT) has been proposed recently \citep{Messer2013}, extending results of \citet{Steinebach1995}. The underlying statistical test investigates the null hypothesis of constant rate. 

Here we investigate the respective filtered derivative process under the alternative of change points in the rate, assuming that also the variance may change simultaneously. We derive the approximation
\begin{align*}
G \sim  \Delta \cdot\left(\Lambda + L\right),
\end{align*}
where the notation '$\cdot$' denotes pointwise multiplication. 
The function $\Lambda$ accounts for the systematic deviation of $G$ in the neighborhood of a change point (section \ref{sect:filtered}). Interestingly, in contrast to similar approaches \citep{Bertrand2000} this function takes the form of a shark's fin here  because both the rate and the variance can change. Thus, we term $\Lambda$ the shark fin function. The height of the shark's fin is proportional to  a scaled rate difference and grows with the bandwidth of $G$. In practical application, the estimators of the point process parameters are not consistent in the neighborhood of a change point. In section \ref{sect:distortion}, we therefore derive the function $\Delta$ termed here the distortion function. It accounts for the lack in consistency and describes the local parameter estimating process  relative to the true scaling of the filtered derivative process.

Note that for convenience all results in the present article are shown here for processes with independent and identically distributed life times but extend directly to a larger class of renewal processes with a certain degree of variability in the variance (renewal processes with varying variance, RPVVs, compare \citet{Messer2014}) using the same proofs.

%%%%%%%%%%%%%%%%%%
\section{The Filtered Derivative Process}\label{sect:filtered}
\subsection{Notation and behavior under the null hypothesis}
The main goal of the MFT proposed in \citet{Messer2013} is to test the null hypothesis $H_0$ of constant rate against the alternative that a process is a piecewise renewal process with a non-empty set of change points $C=\{c_1,\ldots,c_k\}$, and to estimate the change points in case of rejection of the null hypothesis. In this paper we study the behavior of the filtered derivative process under the alternative. To that end, we first define the filtered derivative process and recall a convergence result under $H_0$, which can be used for the statistical test.

Throughout the article we use the following notation: We write a point process $\Phi$ on the positive line as an increasing sequence of events
$0 <  S_1 <  S_2 < S_3 <\cdots$, or alternatively, by its life times $\xi_j:=S_j - S_{j-1}$, $j=2,3,\ldots$, setting $\xi_1 = S_1$, or by the counting process $(N_t)_{t\ge0}$, where
\begin{align}\label{def_nt}
 N_t=\max\{j\ge 1\,|\, S_j\le t\}, \quad t\ge 0,
\end{align}
with the convention $\max \emptyset := 0$. The process $\Phi$ is called a renewal process with square integrable life times (RP) if the associated life times $\{\xi_j\}_{j\ge 1}$ build a sequence of positive, independent and identically distributed (i.i.d.) and square-integrable random variables with $\sigma^2:=\mathbb{V}\!ar(\xi_1)>0$. 
For an RP $\Phi$ with $\mu:=\mathbb{E}[\xi_1]$ and $\sigma^2=\mathbb{V}\!ar[\xi_1]$ we write $\Phi=\Phi(\mu,\sigma^2)$.
The inverse mean $\mu^{-1}$ is termed the rate of $\Phi$.\\ 
%Note that all results derived in this article apply directly to a more general class or renewal processes (renewal processes with varying variance, RPVVs, defined in \citet{Messer2013}), using the same proofs.\\ 
 For $T>0$ let $(\Phi^{(n)})_{n\ge1}:=\Phi |_{(0,nT]}$, where $\Phi |_{(a,b]}$ denotes the restriction of $\Phi$ to the interval $(a,b]$. The value $n$ is required for asymptotic statements throughout this work, which are deduced by letting $n\to\infty$. Here, the total time $nT$ and the location of the change point $nc$ grow linearly in $n$. Let $(N_t^{(n)})_{t\ge0}$ denote the counting process corresponding to $\Phi^{(n)}$. For $T>0$ let $h\in(0,T/2]$ denote a window size  and $\tau_h := [h,T-h]$ an analysis region. 

\begin{defi}\label{def_dht}
Let $\Phi(\mu,\sigma^2)$ be an RP. For $t\in\tau_h$ the filtered derivative process $D^{(n)}:=\left(D_t^{(n)} \right)_{t\in \tau_h}$ is defined as
\begin{align}\label{d_ht}	
D_t^{(n)}:=D_{h,t}^{(n)}:= \frac{\left(N_{n(t+h)}^{(n)}-N_{nt}^{(n)}\right)-\left(N_{nt}^{(n)}-N_{n(t-h)}^{(n)}\right)}{s_{t}^{(n)}},
\end{align}	
where $s_{t}^{(n)}:=s_{h,t}^{(n)}:=\sqrt{2nh\sigma^2/\mu^3}$. 
\end{defi}
Thus, $D_t^{(n)}$ compares the number of events in a left window, $N_{nt}^{(n)}-N_{n(t-h)}^{(n)}$, to the number of events in a right window, $N_{n(t+h)}^{(n)}-N_{nt}^{(n)}$ (Figure \ref{behavior_rate_z}). The process $D^{(n)}$ can indicate changes in the rate because its expectation asymptotically vanishes under the null hypothesis, while systematic deviations from zero are expected when a rate change occurs. More precisely, under the null hypothesis the following weak process convergence result for $D^{(n)}$ was shown in \citet{Steinebach1995} and \citet{Messer2013} for renewal processes and certain generalizations with respect to variability in the variance. Let $D[h,T-h]$ denote the set of all c\`{a}dl\`{a}g (continue \`{a} droite, limite \`{a} gauche) functions on $[h,T-h]$. Further, let $d_{SK}$ denote the Skorokhod metric on $D[h,T-h]$.\\ 
The following result describes the limit behavior of $D$ when no change in the rate occurs. 

\begin{prop}\label{hauptaussage_alternative_no_cp}
Let $\Phi(\mu,\sigma^2)$ be an RP such that $\Phi^{(n)}=\Phi |_{(0,nT]}$. Let $(W_t)_{t\ge 0}$ be a standard Brownian motion. Then it holds in $(D[h,T-h],d_{SK})$ as $n\to\infty$ 
\begin{align}\label{hauptaussage_alternative_no_cp_eq}
 \left(D_t^{(n)}\right)_{t\in\tau_h} \; \stackrel{d}{\longrightarrow} \; \left(\frac{(W_{t+h}-W_t)-(W_t-W_{t-h})}{\sqrt{2h}}\right)_{t\in\tau_h}.
\end{align}
\end{prop}

The expression $\stackrel{d}{\longrightarrow}$ denotes convergence in distribution.
Proposition \ref{hauptaussage_alternative_no_cp} is a special case of Proposition \ref{hauptaussage_alternative} (section \ref{subsec:fdcp}), which describes the behavior of $D^{(n)}$ in the presence of a change point. 
In case of a change point in the rate, $D_t^{(n)}$ systematically deviates from zero in the neighborhood of the change point. Therefore, we introduce an additional centering term in the following subsection in order to obtain convergence in case of a change point.

\subsection{The filtered derivative in case of a change point}\label{subsec:fdcp}
In order to investigate the behavior of $D_t$ under the alternative of change points, we note that a change point at $c$ can only affect $D_t$ within the $h$-neighborhood of $c$, i.e., for  $t \in (c-h,c-h)$. Therefore, investigating one change point extends directly to an  arbitrary number of change points with distances at least $2h$. We thus focus here on the behavior in case of one change point, using the following point process model. The process $\Phi^{(n)}$ starts as the RP $\Phi_1(\mu_1,\sigma_1^2)$ and jumps into $\Phi_{2}(\mu_{2},\sigma_{2}^2)$ at the change point $nc$. 

\begin{figure}[htb]
	\centering
	\includegraphics[scale=0.65]{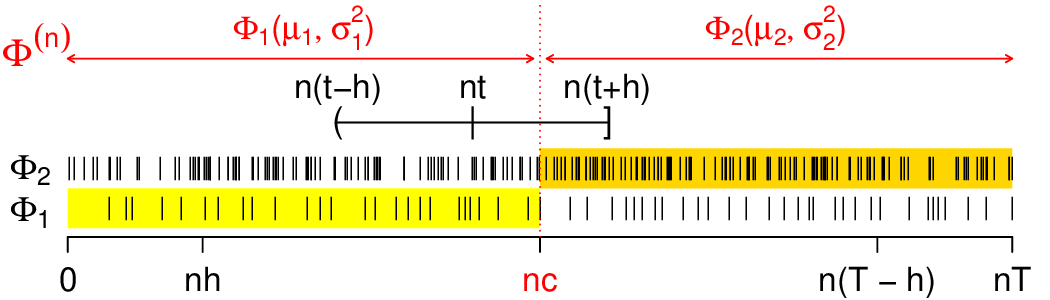}
	\caption{A point process with a change point at $nc$ according to Construction \ref{construction_renewal_model_extension}. Before $nc$, $\Phi^{(n)}$ equals $\Phi_1(\mu_1,\sigma_2^2)$ and after $nc$, it derives from a second RP $\Phi_2(\mu_2,\sigma_2^2)$. The  windows  required for the filtered derivative process at time $nt$ are given by the intervals $(n(t-h),nt]$ and $(nt,n(t+h)]$.} 
	\label{behavior_rate_z}
\end{figure}

\begin{cons}\label{construction_renewal_model_extension}
Let  $c\in(0,T)$ and $n=1,2,\ldots$ 
\noindent Let $\Phi_1(\mu_1,\sigma_1^2)$ and $\Phi_{2}(\mu_{2},\sigma_{2}^2)$ be two independent RPs and set
\begin{align}\label{gamma_model_extended}
\Phi^{(n)}:= \Phi^{(n)}(c) := \Phi_1 |_{(0,nc]}\cup \Phi_2 |_{(nc,nT]}.
\end{align}
The resulting sequence of interest is given as 
$\left(\Phi^{(n)}\right)_{n\ge 1}$ (cmp.~Figure \ref{behavior_rate_z}).
\end{cons}

Thus, a change in the rate occurs if and only if $\mu_1 \neq \mu_{2}$.
In this case of one change point, $D$ will systematically deviate from zero in the $h$-neighborhood of $c$ \citep[cmp.][]{Bertrand2000}. Therefore, we require an additional centering term $m_t$ for process convergence, and an extension of the scaling process $s_t$ as follows

\begin{figure}[htb]
  	\centering
            \includegraphics[scale=0.65]{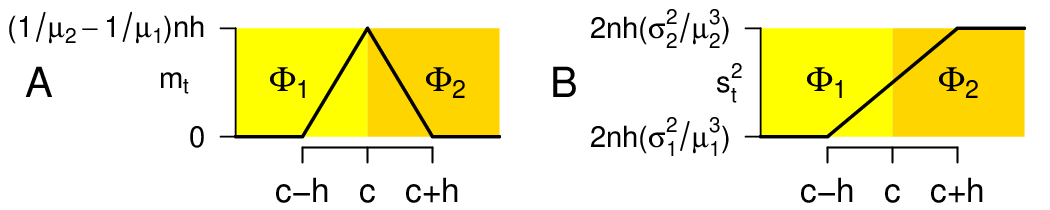}
             \includegraphics[scale=0.65]{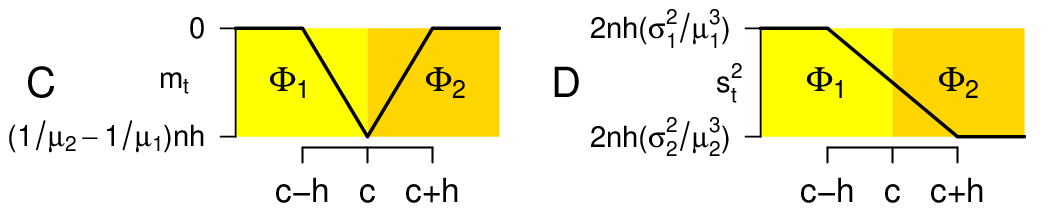}
          \caption{Representation of the expectation function $m_t$ (A, C) and the variance function $s_t^2$ (B, D) in case of a change point at $c$, according to Definition \ref{rescaledfdp}. The expectation $m_t$ vanishes outside $[c-h,c+h]$ and takes its extreme, $m_{c}= (1/\mu_2-1/\mu_1)nh$, at $c$. The function $m_t$ is non-negative if the rate increases (A) and non-positive if the rate decreases (C). The variance function $s_t^2$ equals $2nh(\sigma_1^2/\mu_1^3)$ for $t<c-h$ and $2nh(\sigma_2^2/\mu_2^3)$ for $t>c+h$ and is linearly interpolated in  $[c-h,c+h]$ (B, D). Superscipts $(n)$ are omitted here for convenience in the notation of $m$ and $s$.}
        \label{center_scale_cp}
    \end{figure}

\begin{defi}\label{rescaledfdp}
Let the rescaled filtered derivative process $\Gamma^{(n)}:=\left(\Gamma_t^{(n)}\right)_{t\in\tau_h}$ be defined as
\begin{align}\label{gamma_tilde}	
\Gamma_{t}^{(n)}:=\Gamma_{h,t}^{(n)}:= \frac{\left[(N_{n(t+h)}^{(n)}-N_{nt}^{(n)})-(N_{nt}^{(n)}-N_{n(t-h)}^{(n)})\right]-m_{t}^{(n)}}{s_{t}^{(n)}},
\end{align}	
while for $t\in \tau_h$ the  expectation function $m^{(n)}:=\left(m_t^{(n)}\right)_{t\in\tau_h}$ is zero for $|t-c|>h$ and equals
\begin{equation}\label{true_center}
m_{t}^{(n)} :=m_{h,t}^{(n)}(c) :=
		 n \left(1/\mu_2-1/\mu_1\right) (h-|t-c|)\quad \text{for }\; |t-c|\le h\quad \text{(see Figure \ref{center_scale_cp} A, C)}.
\end{equation}
The variance $(s^{(n)})^2:=\left((s_t^{(n)})^2\right)_{t\in\tau_h}$ is given by $2nh\sigma_1^2/\mu_1^3$ for $t<c-h$, by $2nh\sigma_2^2/\mu_2^3$ for $t>c+h$, and by a linear interpolation (see Figure \ref{center_scale_cp} B, D)
\begin{equation}\label{true_scale}
(s_{t}^{(n)})^2 :=(s_{h,t}^{(n)})^2 :=
		 n \left((t+h-c)\sigma_2^2/\mu_2^3+(c-(t-h))\sigma_1^2/\mu_1^3\right), \quad\text{for }\; |t-c|\le h.
\end{equation}
\end{defi}

Intuitively, the linear interpolation results from the linear shift of the window across time: Assume for example a rate increase (Figure \ref{center_scale_cp} A). If the window is shifted to the right in the interval $(c-h,c)$, only its right half is expected to contain more events. The fraction of the right half for which this is the case increases linearly up to time $c$.  Analogously, the decrease is linear in the interval $(c,c-h)$. For the variance a similar argument holds due to additivity of the variances under independence of the life times.

Similar to the process $D^{(n)}$, also the process $\Gamma^{(n)}$ can be shown to converge weakly in Skorokhod topology to a limit process $L$ in the general setting of a change point, as stated in the following proposition.

\begin{prop}\label{hauptaussage_alternative}
Let $\Phi_1(\mu_1,\sigma_1^2)$ and $\Phi_2(\mu_2,\sigma_2^2)$ be independent RPs with $\mu_1\not=\mu_2$. Let the sequence $(\Phi^{(n)})_{n\ge1}$ result from $\Phi_1$ and $\Phi_2$ according to Construction \ref{construction_renewal_model_extension}, and let $\Gamma^{(n)}$ be the associated rescaled filtered derivative process.  Let $(W_t)_{t\ge 0}$ be a standard Brownian motion, and the limit process $L:=(L_{t})_{t\in\tau_h}$ be given as
\begin{align}\label{limit_process_cp}
L_{t}:=L_{h,t}(c)	&:= \begin{cases}
				\frac{(W_{t+h}-W_t)-(W_t-W_{t-h})}{\sqrt{2h}}, 	& \textrm{if }\; |t-c|>h,\vspace{1em}\\
				\frac{\sqrt{\sigma_2^2/\mu_2^3}(W_{t+h}-W_c) + \sqrt{\sigma_1^2/\mu_1^3}[(W_c-W_t)-(W_t-W_{t-h})]}{s_{t}^{(1)}}, & \textrm{if }\; c-h\le t\le c,\vspace{1em}\\
				\frac{\sqrt{\sigma_2^2/\mu_2^3}[(W_{t+h}-W_t)-(W_t-W_{c})] - \sqrt{\sigma_1^2/\mu_1^3}(W_{c}-W_{t-h})}{s_{t}^{(1)}}, & \textrm{if }\; c<t\le c+h.
				\end{cases}	
\end{align}
Then it holds in $(D[h,T-h],d_{SK})$ as $n\to\infty$ 
\begin{align*}
\Gamma^{(n)} \; \stackrel{d}{\longrightarrow} \; L.
\end{align*}
\end{prop}

 Elementary calculations show that the marginals $L_t$ are standard normally distributed. Note that Proposition \ref{hauptaussage_alternative_no_cp} describes the special case where $\mu_1=\mu_2$ and $\sigma_1^2=\sigma_2^2$, because for all $t\in\tau_h$, we obtain $m_{t}^{(n)}=0$, $s_{t}^{(n)}=(2nh\sigma_1^2/\mu_1^3)^{1/2}$, $L_{t}=[(W_{t+h}-W_t)-(W_t-W_{t-h})]/(2h)^{1/2}$ and $\Gamma_{t}^{(n)}$ equal to the left hand side in equation (\ref{hauptaussage_alternative_no_cp_eq}). 
The idea of the proof of Proposition \ref{hauptaussage_alternative} is similar, but it
relies on joint process convergence of the rescaled counting processes associated with $\Phi_1$ and $\Phi_2$, see Appendix \ref{subsec:proof_alt}.

\subsection{The Shark Fin Function}

Proposition \ref{hauptaussage_alternative} states that asymptotically the following equality in distribution holds
\begin{align}\label{lambda}
 D^{(n)} \; \sim \; \Lambda^{(n)}+ L, \quad \text{with} \quad \Lambda^{(n)}:= m^{(n)}/s^{(n)}.
\end{align}
In order to understand the process $D^{(n)}$ we investigate $\Lambda^{(n)}$. In case of a rate change the expectation function $m^{(n)}$ has the shape of a hat (Figure \ref{center_scale_cp} A). If the variance changes additionally, the function $\Lambda^{(n)}$ resembles a shark's fin and is therefore termed here the shark fin function. 

\begin{figure}[htb]
	\centering
	\includegraphics[scale=0.65]{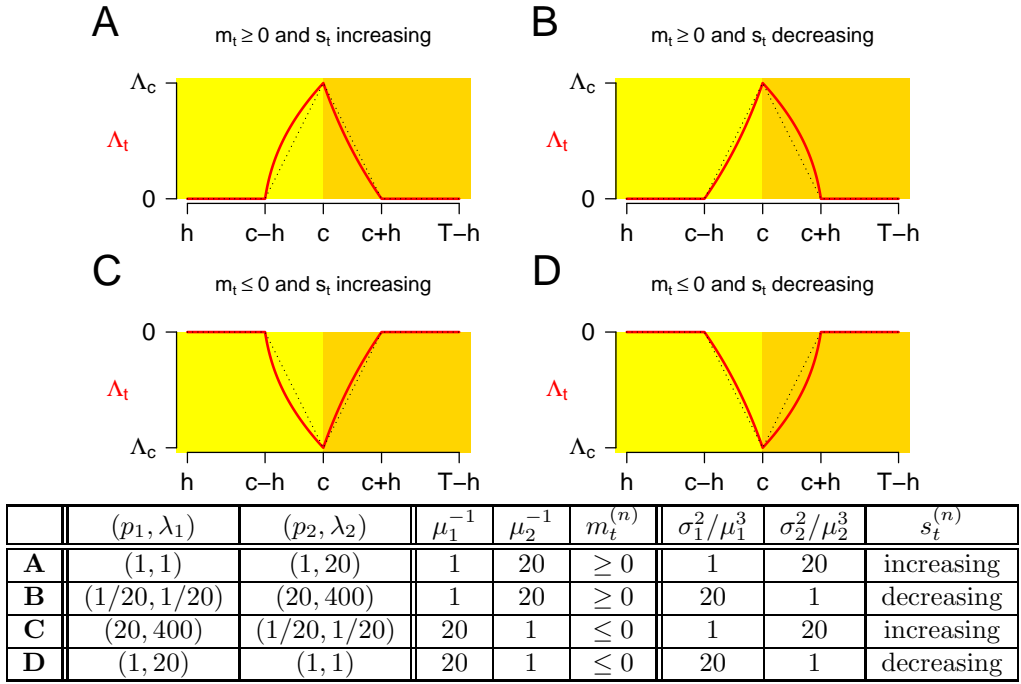}
	\vspace{0.3em}
	\begin{small}
		\begin{tabular}{| c || c | c || c | c | c || c | c | c |}
			\hline
			&$(p_1,\lambda_1)$	& $(p_2,\lambda_2)$ & $\mu_1^{-1}$ & $\mu_2^{-1}$ & $m_{t}^{(n)}$ & $\sigma_1^2/\mu_1^3$ & $\sigma_2^2/\mu_2^3$ & $s_{t}^{(n)}$ \\
			\hline
			\hline
			\textbf{A} 	& 	$(1,1)$			&   $(1,20)$    		& $1$	&		$20$	    &		$\ge0$	  & $1$				&	$20$				   & increasing	\\
			\hline
			\textbf{B} 	& 	$(1/20,1/20)$		&   $(20,400)$    	& $1$	&		$20$	    &		$\ge0$	  & $20$				&	$1$				   & decreasing	\\
			\hline
			\textbf{C}	& 	$(20,400)$		&   $(1/20,1/20)$    	& $20$	&		$1$	    &		$\le0$	  & $1$				&	$20$				   & increasing	\\
			\hline
			\textbf{D} & 	$(1,20)$			&   $(1,1)$    		& $20$	&		$1$	    &		$\le0$	  & $20$				&	$1$				   & decreasing	\\
			\hline
		\end{tabular}
	\end{small}
	\caption{Analysis of the shark fin function $\Lambda_{t}$ (solid), for the case of a change point at  $c$. The dotted line marks the scaled hat function $m_{t}/s_{c}$. The shape of the shark fin function depends on the structure of the expectation function $m_{t}$ and the standard deviation function $s_{t}$. For $m_{t}\ge0$ and $s_{t}^2$ increasing, the shark is going west (A).  For $m_{t}\ge0$ and $s_{t}^2$ decreasing, the shark is heading east (B). For $m_{t}\le0$, the shark is swimming upside down and is oriented towards the same directions (C and D). The expectation and standard deviation functions refer to point processes whose life times before the change point are i.i.d. $\Gamma(p_1,\lambda_1)$ distributed and those after the change point are i.i.d. $\Gamma(p_2,\lambda_2)$ distributed. The parameters are given in the upper table. Further parameters are $T=1000$, $c=500$, $h=150$ and $n=1$. Superscripts ${(n)}$ are omitted  for convenience.}
	\label{shark}
\end{figure}

We show examples of such shark fin functions in Figure \ref{shark} and give a proof in Lemma \ref{lemm_shark}.
Equation (\ref{shark_max}) states that the shark fin function takes its largest deviation from zero at time $c$. If $m^{(n)}\ge 0$ and $s^{(n)}$ increasing, the shark is heading west (Figure \ref{shark} A, equation (\ref{shark_a})), whereas in case of (\ref{shark_b}), the shark is heading east (Figure \ref{shark} B). For $m^{(n)}\le 0$ analogous relations hold, the shark is heading in the same directions, but turned upside down (Figure \ref{shark} C and D). Note also that if the standard deviation $s^{(n)}$ is constant over time, the shark fin function $\Lambda^{(n)}$ has a hat shape, i.e., is piecewise linear.

\begin{lemm}\label{lemm_shark}
For $c\in(0,T)$ and  $t\in\tau_h$ let $m^{(n)}$ and $s^{(n)}$ be as in Definition \ref{rescaledfdp} and $\Lambda^{(n)}= m^{(n)}/s^{(n)}$.  
Then $\Lambda^{(n)}$ is a continuous function with $\Lambda^{(n)}=0$
 for $t \notin (c-h,c+h]$.  If $\mu_1=\mu_2$ it also is $\Lambda^{(n)}=0$ for $t \in (c-h,c+h]$. If $\mu_1\not=\mu_2$ we separate four cases  for $t \in (c-h,c+h]$: 
For $m^{(n)}\ge0$ and $s^{(n)}$ increasing (Figure \ref{shark} A), 
\begin{align}\label{shark_a}
\Lambda_{t}^{(n)}& \quad\textrm{is}\quad
\begin{cases}
\textrm{concave and strictly increasing for t} \in [c-h,c],\\
\textrm{convex and strictly decreasing for t} \in (c,c+h].
\end{cases}
\end{align}
For $m^{(n)}\ge0$ and $s^{(n)}$ decreasing (Figure \ref{shark} B), 
\begin{align}\label{shark_b}
\Lambda_{t}^{(n)}& \quad\textrm{is}\quad
\begin{cases}
\textrm{convex and strictly increasing for t} \in [c-h,c],\\
\textrm{concave and strictly decreasing for t} \in (c,c+h].
\end{cases}
\end{align}
For $m^{(n)}\le0$, expressions (\ref{shark_a}) and (\ref{shark_b}) hold true, but with 'convex' and 'concave' as well as 'increasing' and 'decreasing' exchanged.\\
Further, because $m_t^{(n)}$ is of order $nh$ and $s_t^{(n)}$ is of order $(nh)^{1/2}$ for $|t-c|<h$, we find that $\Lambda_t^{(n)}$ is of order $(nh)^{1/2}$ for $|t-c|<h$.
\end{lemm}

\textbf{Proof of Lemma \ref{lemm_shark}}: 
Continuity is clear because both the numerator and the denominator are continuous. 
For $t \notin (c-h,c+h]$ it is $m^{(n)}=0$ such that $\Lambda^{(n)}=0$. The same holds for $t \in (c-h,c+h]$ when $\mu_1=\mu_2$. For $t \in (c-h,c+h]$ with $\mu_1\not=\mu_2$ we deduce the case $m^{(n)}\ge 0$ and $s^{(n)}$ increasing. For $t\in (c-h,c]$ both functions $m^{(n)}$ and $s^{(n)}$ are strictly increasing in $t$. While $m_t$ is of order $t$, $s_t$ is of order $t^{1/2}$, see equations (\ref{true_center}) and (\ref{true_scale}). Thus, the shark fin function $\Lambda^{(n)}$ is strictly increasing and of order $t^{1/2}$, and therefore describes a concave function for $t\in(c-h,c]$. For $t\in (c,c+h]$, $m^{(n)}$ is strictly decreasing and of order $t$, so that $\Lambda^{(n)}$ is strictly decreasing of order $t^{1/2}$, which describes a convex function. The other cases follow by similar arguments. \hfill $\Box$\\

Note that, if $\mu_1\not=\mu_2$ because $\Lambda_t^{(n)}$ is defined for $t \in \tau_h$ we find
\begin{align}\label{shark_max}
\arg\max_t\left| \Lambda_{t}^{(n)}\right| = c
\end{align}	
for $c\in\tau_h$, $\arg\max |\Lambda_{t}^{(n)}|=h$ if $c\in (0,h)$ and $\arg\max |\Lambda_{t}^{(n)}|=T-h$ if $c\in (T-h,T)$. Note further that Lemma \ref{lemm_shark} can be generalized to multiple change points with distance  at least $2h$, in which case $\Lambda^{(n)}$ describes multiple, successive shark fin functions. %If change points lie closer together, more technical effort is necessary in order to formulate an appropriate convergence result. 

\paragraph{Detection Probability in Change Point Estimation}
The fact that $\Lambda^{(n)}$ takes its maximal deviation from zero at the change point $c$ can be used for change point estimation and for a rough evaluation of the detection probability of a change point.
In practice, the null hypothesis of constant rate is rejected if the filtered derivative $D^{(1)}$ exceeds a threshold $Q$, which can be derived  by Monte Carlo simulation, compare e.g.~\citet{Messer2013}. If the null hypothesis is rejected, an estimate of a change point $c$ is given as $\hat c:=\arg\max_{t\in\tau_h} |D^{(1)}|$. For multiple change points, successive argmax-type estimation methods are applied  \citep[cmp. ][]{Carlstein1988,Duembgen1991,Antoch1994,Antoch1997,Bertrand2000,Bertrand2011,Messer2013,Kirch2014}. 

The construction $D^{(n)} =\Lambda^{(n)}+\Gamma^{(n)}$ gives a simple bound for the detection probability of a change point $c\in\tau_h$. According to Proposition \ref{hauptaussage_alternative} and equations (\ref{true_center}) and (\ref{true_scale}), we find asymptotically 
\begin{align}\label{bound}
D_{c}^{(n)}\sim\Lambda_{c}^{(n)}+L_{c}\sim N\left(\frac{1/\mu_2-1/\mu_1}{(\sigma_2^2/\mu_2^3+\sigma_1^2/\mu_1^3)^{1/2}}\, (nh)^{1/2},1\right).
\end{align}
For rate increases $\mu_2^{-1}>\mu_1^{-1}$, we find $m_c>0$ and $D_c>0$, such that $P(\max_{t\in\tau_h}|D_{t}^{(n)}|>Q)\ge P(D_c^{(n)}>Q)$. Analogous results apply for rate decreases. This implies asymptotically
\begin{align}\label{bound2}
P\left(\max_{t\in\tau_h}|D_{t}^{(n)}|>Q\right) 	& \ge  1 - F \left(Q-\frac{|1/\mu_2-1/\mu_1|}{(\sigma_2^2/\mu_2^3+\sigma_1^2/\mu_1^3)^{1/2}}\, (nh)^{1/2}\right),
\end{align}
where $F$ denotes the  distribution function of the standard normal distribution. 

Note that the right hand side of equation (\ref{bound}) implies that the height of the shark is proportional to the scaled rate differences and grows with the bandwidth of $D$.

%In order to show show inequality (\ref{bound}), we differentiate two cases. For an increase in the rate $\mu_2^{-1}>\mu_1^{-1}$, we bound the detection probability as
%\begin{align}\label{proof_bound_1}
%P\left(\max_{t\in\tau_h}|D_{h,t}^{(n)}|>Q\right) 	& \ge P\left(D_{h,c}^{(n)}>Q\right)
%									  \approx 1 - F \left(Q-\frac{1/\mu_2-1/\mu_1}{(\sigma_2^2/\mu_2^3+\sigma_1^2/\mu_1^3)^{\nicefrac{1}{2}}}\, (nh)^{\nicefrac{1}{2}}\right).
%\end{align}
%Analogously, for a decrease in the rate, i.e., $\mu_2^{-1}<\mu_1^{-1}$, we bound
%\begin{align}\label{proof_bound_2}
%P\left(\max_{t\in\tau_h}|D_{h,t}^{(n)}|>Q\right) 	 \ge P\left(D_{h,c}^{(n)}<-Q\right)
%									&\approx F \left(-Q-\frac{1/\mu_2-1/\mu_1}{(\sigma_2^2/\mu_2^3+\sigma_1^2/\mu_1^3)^{\nicefrac{1}{2}}}\, (nh)^{\nicefrac{1}{2}}\right)\nonumber\\
%									& = 1- F\left(Q-\frac{1/\mu_1-1/\mu_2}{(\sigma_2^2/\mu_2^3+\sigma_1^2/\mu_1^3)^{\nicefrac{1}{2}}}\, (nh)^{\nicefrac{1}{2}}\right).
%\end{align}
%Both inequalities (\ref{proof_bound_1}) and (\ref{proof_bound_2}) yield (\ref{bound}).

\section{The Distortion -- Estimation of Process Parameters}\label{sect:distortion}

The definition of the filtered derivative process $D^{(n)}$ as in equation (\ref{d_ht}) relies on the assumption that the theoretical standard deviation $s^{(n)}$ is known. However, $s^{(n)}$ depends on the  point process parameters $\mu_1, \mu_2$, $\sigma_1^2$ and $\sigma_2^2$, which typically need to be estimated in practical application. 
Note that the filtered derivative is a local statistic, such that $s^{(n)}$ itself is also a time dependent function in case of rate changes, see  definition (\ref{true_scale}).
We discuss the behavior of the filtered derivative process when replacing $s^{(n)}$ by a time dependent estimator $\hat s^{(n)}$ proposed in \citet{Messer2013}.
There, consistency was shown under $H_0$. 
Here, we deduce the asymptotics of the process $\hat s^{(n)}$ under $H_A$. 
The estimator is not consistent, but deviates from the true scaling $s^{(n)}$ in the $h$-neighborhood of a change point. However, both functionals $\hat s^{(n)}$ and $s^{(n)}$ are of the same magnitude and their asymptotic relation is termed here the distortion $\Delta$. The latter can be interpreted as the amount of error that results from a bias in the parameter estimation close to a change point. For a similar phenomenon in the setting of sequences of random variables compare \citet{Kirch2014}.

 For all $t\in\tau_h$, the estimator $\hat s^{(n)}$ is given by
\begin{align}\label{schaetzer_s}
\left(\hat s_t^{(n)}\right)^2 := \left(\hat s_{h,t}^{(n)}\right)^2 := \left(\frac{\hat\sigma_{ri}^2(nh,nt)}{\hat\mu_{ri}^3(nh,nt)} + \frac{\hat\sigma_{le}^2(nh,nt)}{\hat\mu_{le}^3(nh,nt)} \right)nh,
\end{align}
where $\hat\mu_{ri}(nh,nt)$ and $\hat\sigma_{ri}^2(nh,nt)$ (or $\hat\mu_{le}(nh,nt)$ and $\hat\sigma_{le}^2(nh,nt)$) denote the empirical mean and variance of all life times whose corresponding point events lie in the right window $(nt,n(t+h)]$ (or the left window $(n(t-h),nt]$, respectively). If no life times can be found in the respective intervals,  the estimators are set to zero.

Replacing  $s_{t}^{(n)}$ with this estimator $\hat s_t^{(n)}$, we study the convergence of a new process defined as
\begin{align}\label{g_ht}	
G_{t}^{(n)}:= \frac{(N_{n(t+h)}^{(n)}-N_{nt}^{(n)})-(N_{nt}^{(n)}-N_{n(t-h)}^{(n)})}{\hat s_t^{(n)}}.%\quad \text{for} \quad t\in\tau_h.
\end{align}	
Under the null hypothesis of no change point (i.e., $\mu_1=\mu_2$), the following convergence result is provided in \citet{Messer2013}. 
\begin{prop}\label{hauptaussage_alternative_no_cp_ght}
Let $\Phi_1(\mu_1,\sigma_1^2)$ be an RP, such that $\Phi^{(n)}=\Phi_1 |_{(0,nT]}$. 
Then, we have in $(D[h,T-h],d_{SK})$ as $n\to\infty$ 
\begin{align}\label{hauptaussage_alternative_no_cp_eq_ght}
G^{(n)} \stackrel{d}{\longrightarrow} \; \left(\frac{(W_{t+h}-W_t)-(W_t-W_{t-h})}{\sqrt{2h}}\right)_{t\in\tau_h}.
\end{align}
\end{prop}

The proof relies on the strong consistency of the estimator $\hat s^{(n)}$ under the null hypothesis, i.e., that it holds uniformly almost surely $ s^{(n)}/\hat s^{(n)}\to 1$ as $n\to\infty$.\\ 

In the general case of a change point the relation  $s^{(n)}/\hat s^{(n)}$ does not converge to unity, but to a deterministic function $\Delta:=s^{(1)}/\hat s^{(1)}$.
For all $t\in\tau_h$ let $\tilde s^{(n)}$ be
\begin{align}\label{schaetzer_stilde}
\left(\tilde s_t^{(n)}\right)^2 := \left(\tilde s_{h,t}^{(n)}\right)^2 := \left(\frac{\sigma_{ri}^2(h,t)}{\mu_{ri}^3(h,t)} + \frac{\sigma_{le}^2(h,t)}{\mu_{le}^3(h,t)} \right)nh,
\end{align}
with $\mu_{ri}(h,t)=\mu_1$ for $t\le c-h$, $\mu_{ri}(h,t)=\mu_2$ for $t>c$, and 
\begin{equation}\label{mu_ri_cp}
\mu_{ri}(h,t)=h\mu_1\mu_2/((c-t)\mu_2+(t+h-c)\mu_1) \quad \text{for}\quad t\in(c-h,c],
\end{equation}
and analogously for $\mu_{le}$.
For $\sigma_{ri}$ we set $\sigma_{ri}^2(h,t)=\sigma_1^2$ for $t\le c-h$, $\sigma_{ri}^2(h,t)=\sigma_2^2$ for $t>c$ and
\begin{equation}\label{sig_ri_cp}
\sigma^2_{ri}(h,t)=\frac{\mu_1\mu_2(t+h-c)(c-t)[(\sigma_1-\sigma_2)^2+(\mu_1+\mu_2)^2]+[(t+h-c)\mu_1\sigma_2+(c-t)\mu_2\sigma_1]^2}{[(c-t)\mu_2+(t+h-c)\mu_1]^2} 
\end{equation}
for $t\in(c-h,c]$, and analogously for $\sigma^2_{le}$.
Let $d_{\|\cdot\|}$ denote the supremum norm.
The following Lemma states that $\Delta$ describes the asymptotic error induced by the estimator $\hat s^{(n)}$. 

\begin{lemm}\label{conv_delta}
	Let $\mu_1, \mu_2, \sigma_1^2, \sigma_2^2>0$.  Let $\Phi_1(\mu_1,\sigma_1^2)$ and $\Phi_2(\mu_2,\sigma_2^2)$ be independent RPs and $c\in (0,T)$, so that the sequence $(\Phi^{(n)})_{n\ge1}$ results from $\Phi_1$ and $\Phi_2$ according to Construction \ref{construction_renewal_model_extension}. 
	Let $s^{(n)}$, $\hat s^{(n)}$ and $\tilde s^{(n)}$ 
 as defined in (\ref{true_scale}), (\ref{schaetzer_s}) and  (\ref{schaetzer_stilde}). Then we have in $(D[h,T-h],d_{\|\cdot\|})$ almost surely as $n\to\infty$ 
\begin{align}\label{conv_distortion_ha}
\left(\frac{s_t^{(n)}}{\hat s_t^{(n)}}\right)_{t\in\tau_h} \longrightarrow \;\; \left(\frac{s_t^{(1)}}{\tilde s_t^{(1)}}\right)_{t\in\tau_h} = \left(\Delta_{t}\right)_{t\in\tau_h}.
\end{align}
\end{lemm}	

\begin{figure}[htbp]
  	\begin{center}
	 \includegraphics[scale=0.6]{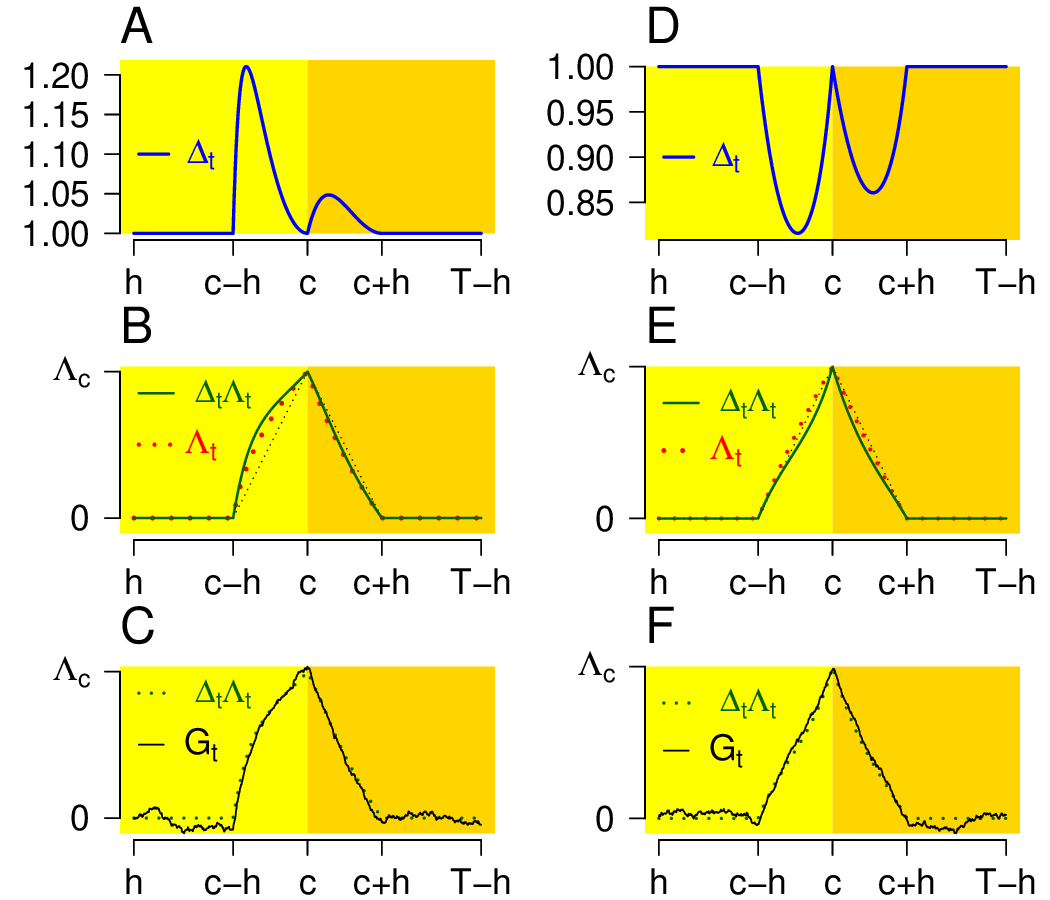}		
	\end{center}
	 \caption{Two examples of the function $G_{t}$ and  its connection to the shark fin function $\Lambda_t$ and the distortion $\Delta_t$. The underlying point process on $[0,1000]$ starts in a process with i.i.d.~$\Gamma(p_1,\lambda_1)$-distributed life times and jumps into a process with i.i.d.~$\Gamma(p_2,\lambda_2)$-distributed life times at time $c=500$. Panels A,D: The distortion $\Delta_{t}$. B,E: The distorted shark fin function  $\Delta_{t}\Lambda_{t}$ (solid), the undistorted shark fin function $\Lambda_t$ (dotted, thick) and the hat function (dotted, thin). C,F: The process $G_{t}$ (solid) that fluctuates around the distorted shark fin function (dotted). For panels A-C, $(p_1,\lambda_1)=(1,5), (p_2,\lambda_2)=(1/4,5)$, resulting in $(\mu_1,\sigma_1^2)=(1/5,1/25)$ and $(\mu_2,\sigma_2^2)=(1/20,1/100)$. For panels D-F, $(p_1,\lambda_1)=(2,10), (p_2,\lambda_2)=(2,20)$, resulting in $(\mu_1,\sigma_1^2)=(1/5,1/50)$ and $(\mu_2,\sigma_2^2)=(1/20,1/200)$. The window size was $h=150$. }
	\label{hat_shark_g}	
    \end{figure}

The proof is given in Appendix \ref{proof:conv_delta}.
Note that this Lemma states that the estimator $\hat s^{(n)}$ asymptotically equals $\tilde s^{(n)}$ almost surely in $(D[h,T-h],d_{\|\cdot\|})$. The distortion $\Delta$ is continuous and depends on the  process parameters $\mu_1,\mu_2,\sigma_1^2$ and $\sigma_2^2$ (see Figure \ref{hat_shark_g} A,D for examples). 

Considering the distortion term for  applications in which the process parameters need to be estimated, we find the following convergence of the filtered derivative process $G^{(n)}$. 

\begin{prop}\label{hauptaussage_alternative_ght}
Let $\mu_1, \mu_2, \sigma_1^2, \sigma_2^2>0$.  Let $\Phi_1(\mu_1,\sigma_1^2)$ and $\Phi_2(\mu_2,\sigma_2^2)$ be independent RPs and $c\in (0,T)$ be a change point, so that the sequence $(\Phi^{(n)})_{n\ge1}$ results from $\Phi_1$ and $\Phi_2$ according to Construction \ref{construction_renewal_model_extension}. 
Then, for $G^{(n)}$, $\Delta$, $\Lambda^{(n)}$ and $L$ as defined in (\ref{g_ht}), (\ref{lambda}), (\ref{conv_distortion_ha}) and (\ref{limit_process_cp}), we have in $(D[h,T-h],d_{SK})$ as $n\to\infty$ 
\begin{align*}
 G^{(n)} - \Delta \Lambda^{(n)}\; \stackrel{d}{\longrightarrow} \; \Delta L.
\end{align*}
\end{prop}

\noindent\textbf{Proof of Proposition \ref{hauptaussage_alternative_ght}}: 
Since $\Delta \Lambda^{(n)}={m^{(n)}}/{\hat s^{(n)}}$, the claim follows directly from
\begin{align*}
\Gamma_t^{(n)}= \left(G_t^{(n)} -\frac{m_t^{(n)}}{\hat s_t^{(n)}}\right) \frac{\hat s_t^{(n)}}{s_t^{(n)}}
\end{align*}
and due to the weak convergence $\Gamma^{(n)}\to L$	 as stated in Proposition \ref{hauptaussage_alternative} and the almost sure convergence $ s^{(n)}/\hat s^{(n)} \to \Delta$ as in Lemma \ref{conv_delta} by applying Slutsky's theorem.\hfill $\Box$\\

As a corollary we note that if $\mu_1=\mu_2$, even if $\sigma_1^2\not=\sigma_2^2$, we find that $\Lambda_t^{(n)}=0$ and $\Delta=1$, which can be obtained by elementary calculations. 

\begin{coro}
	Let  $c\in(0,T)$ and $n=1,2,\ldots$ 
	\noindent Let $\Phi_1(\mu,\sigma_1^2)$ and $\Phi_{2}(\mu,\sigma_{2}^2)$ be two independent RPs with $\sigma_1^2 \neq \sigma_{2}^2$ and set
	$\Phi^{(n)}:= \Phi^{(n)}(c) := \Phi_1 |_{(0,nc]}\cup \Phi_2 |_{(nc,nT]}.$
	Then, for $G^{(n)}$ and $L$ as defined in (\ref{g_ht}) and (\ref{limit_process_cp}), we have in $(D[h,T-h],d_{SK})$ as $n\to\infty$ 
	\begin{align*}
	G^{(n)} \; \stackrel{d}{\longrightarrow} \; L.
	\end{align*}
\end{coro}	

Note that the impact of the distortion function  may theoretically become arbitrarily large for extreme parameter constellations (up to $20$\% of the shark fin function in Figures \ref{hat_shark_g} A,D). However, because the estimators are derived locally and separately in each window half, the estimation at the change point $c$ is consistent 
 and the distortion is unity. As a consequence, the estimation error caused by inconsistent parameter estimation in practical application is typically  small  because the shark fin function takes its largest deviation at $c$.

\section{Summary}

We extend a convergence result of a filtered derivative process described by \citet{Steinebach1995} and \citet{Messer2013} that can be used for change point analysis in point processes. Usually, for purposes of statistical hypothesis testing, the behavior of the filtered derivative process $G$ is analyzed under the null hypothesis. In the present setting it converges weakly to a zero-mean, unit variance Gaussian process $L$ (equation (\ref{limit_process_cp}), upper case), i.e.,
\begin{equation*}
G^{(n)}\sim L.
\end{equation*}
Zero expectation results from a constant rate. Since the parameter estimators are consistent under the null hypothesis, no additional term is required to describe the limit behavior of $G^{(n)}$.\\
 The main purpose of this paper was to describe the  behavior of $G^{(n)}$ under the alternative of one change point. Proposition \ref{hauptaussage_alternative_ght} states that we can approximate (roughly) 
 \begin{align}\label{G_final}
 G^{(n)} \sim \Delta\cdot \left(\Lambda^{(n)}+ L\right).
 \end{align}
The systematic term $\Lambda^{(n)}$ describes the expectation of the filtered derivative, which systematically deviates from zero in the neighborhood of a change point. Interestingly, this deviation does not simply take the form of a hat, but of a shark's fin. This is caused by the assumption that both the rate and the variance may change at a change point.  
 In practice, this shape is distorted further when the process parameters need to be estimated. The distortion function $\Delta$ accounts for the lack in consistency in estimation of point process parameters in the neighborhood of a change point.
 
In summary, the first part in (\ref{G_final}), $\Delta \Lambda^{(n)}$ describes the deterministic, distorted shark fin function (Figure \ref{hat_shark_g}). %, where the distortion is caused by the need to estimate the process parameters and the shark fin function is a systematically deformed hat whose deformation originates from the scaling that is required for process convergence.  
 The second part, $\Delta L$, describes a random fluctuation with zero expectation and variance given as the squared distortion. As a consequence of the local nature of $G^{(n)}$, this result applies automatically to multiple change points separated by at least $2h$. 
Our results also suggest that in practical application the shape of the potentially distorted shark's fin typically neither affects the detection and estimation of change points, nor the lower bound of the detection probability.

\section*{Acknowledgements}
This work was supported by the German Federal Ministry of Education and Research (BMBF) within the framework of the e:Med research and funding concept (grant number: 01ZX1404B) and by the Priority Program 1665 of the DFG. We thank Brooks Ferebee for helpful comments on the manuscript.

\bibliography{MesserSchneider}

\appendix
\section{Appendix}
Unless otherwise specified, we use the following notation (compare Construction \ref{construction_renewal_model_extension}): Let $T>0$, $h\in(0,T/2]$, $t\in\tau_h$ and $c\in(0,T)$.
Further, let $\{\xi_{1,j}\}_{j \ge 1}$, $\{\xi_{2,j}\}_{j \ge 1}$ and $\{\xi_j^{(n)}\}_{j\ge 1}$ denote the sequences of life times that correspond to $\Phi_1$, $\Phi_2$ and to the compound process $\Phi^{(n)}$, respectively. Analogously, let $(N_{1,t})_{t\ge0}$, $(N_{2,t})_{t\ge0}$ and $(N_t^{(n)})_{t\ge0}$ denote the associated counting processes (see equation (\ref{def_nt})).  Further, let $(W_{1,t})_{t\ge 0}$ and $(W_{2,t})_{t\ge 0}$ be independent standard Brownian motions. 

\subsection{Proof of Proposition \ref{hauptaussage_alternative}}\label{subsec:proof_alt}
Outline: We show the joint convergence in distribution of the rescaled counting processes $(N_{1,t})_{t}$ and $(N_{2,t})_{t}$ to a function of  $(W_{1,t})_{t}$ and $(W_{2,t})_{t}$ (compare \ref{konvergenz_z_schlange_gemeinsam}). Then, at time $t$ both processes refer to the information of the entire time interval $(0,t]$. In a second step, the processes are continuously mapped to the scenario of the two windows $(t-h,t]$ and $(t,t+h]$ which refers to the filtered derivative process $(\Gamma_t)_t$. 

\noindent\textbf{Proof of Proposition \ref{hauptaussage_alternative}:}

For $i=1,2$ let the rescaled random walk $(X_{i,t}^{(n)})_{t\ge 0}$ and the rescaled counting process $(Z_{i,t}^{(n)})_{t\ge 0}$ concerning $\Phi_i$ be given as
\begin{align}
X_{i,t}^{(n)}  := \frac{1}{\sigma_i\sqrt{n}}\sum_{j=1}^{[nt]}(\xi_{i,j}-\mu_i)\quad \text{and}\quad
Z_{i,t}^{(n)}  := \frac{N_{i,n t} - n t/\mu_i}{\sqrt{n\sigma_i^2/\mu_i^3}},
\end{align}
for $t\ge 0$.
According to Donsker's theorem (in the case of RPVVs apply \citet[Proposition A.8.]{Messer2013}), we find in $(D[0,\infty),d_{SK})$ as $n\to\infty$ that 
$$(X_{i,t}^{(n)})_{t\ge 0}\stackrel{d}{\longrightarrow} (W_{i,t})_{t\ge 0} \quad \text{for}\quad i=1,2,$$ 
implying weak convergence of $(Z_{i,t}^{(n)})_{t\ge 0}$, i.e., it holds in $(D[0,\infty),d_{SK})$ as $n\to\infty$ that $(Z_{i,t}^{(n)})_{t\ge 0}\stackrel{d}{\longrightarrow} (W_{i,t})_{t\ge 0}$ for $i=1,2$, as stated in \citet[Theorem 14.6.]{Billingsley1999}.

We use a different scaling and set
\begin{align*}
\widetilde Z_{i,t}^{(n)}  := \frac{N_{i,n t} - n t/\mu_i}{s_{t}^{(n)}}, \quad t\ge 0,
\end{align*}
where $s_{t}^{(n)}, t\in [0,\infty)$ is given in Definition \ref{rescaledfdp}. Then for $i=1,2$, we find in $(D[0,\infty),d_{SK})$ for $n\to\infty$
\begin{align*}
\left(\widetilde Z_{i,t}^{(n)} \right)_{t\ge 0}\stackrel{d}{\longrightarrow}  \left(\frac{\sqrt{\sigma_i^2/\mu_i^3}}{s_{t}^{(1)}} W_{i,t} \right)_{t\ge 0}
\end{align*}
because 
$\left({\sqrt{n}\sqrt{\sigma_i^2/\mu_i^3}}/{s_{t}^{(n)}}\right)_{t} = \left({\sqrt{\sigma_i^2/\mu_i^3}}/{s_{t}^{(1)}}\right)_{t} 
$
is continuous in $t$ and does not depend on $n$.

Let now $(\widetilde Z_{1,t}^{(n)})_{t\ge0}$ and $(\widetilde Z_{2,t}^{(n)})_{t\ge0}$ denote the processes derived from $\Phi_1$ and $\Phi_2$, respectively. Due to independence of $\Phi_1$ and $\Phi_2$, we obtain joint convergence in $(D[0,\infty)\times D[0,\infty), d_{SK} \otimes d_{SK})$ for $n\to\infty$
\begin{align}\label{konvergenz_z_schlange_gemeinsam}
\left(\left(\widetilde Z_{1,t}^{(n)} \right)_{t\ge 0}, \left(\widetilde Z_{2,t}^{(n)} \right)_{t\ge 0}\right) \stackrel{d}{\longrightarrow}
\left(\left(\frac{\sqrt{\sigma_1^2/\mu_1^3}}{s_{t}^{(1)}} W_{1,t} \right)_{t\ge 0},\left(\frac{\sqrt{\sigma_2^2/\mu_2^3}}{s_{t}^{(1)}} W_{2,t} \right)_{t\ge 0} \right).
\end{align}
We consider the continuous map $\varphi: (D[0,\infty)\times D[0,\infty), d_{SK} \otimes d_{SK}) \to (D[h,T-h],d_{SK})$ given by
\begin{align*}
((f(t))_{t\ge0},(g(t))_{t\ge0})\\
\stackrel{\varphi}{\longmapsto}   &\left( \begin{array}[c]{l}
(f(t+h)-f(t)) - (f(t)-f(t-h))\mathbbm{1}_{[h,c-h)}(t) 	\\
+(g(t+h)-g(c))+(f(c)-f(t)) - (f(t)-f(t-h))\mathbbm{1}_{[c-h,c)}(t) \\
+(g(t+h)-g(t))-(g(t)-g(c)) - (f(c)-f(t-h))\mathbbm{1}_{[c,c+h)}(t)\\
+(g(t+h)-g(t)) - (g(t)-g(t-h))\mathbbm{1}_{[c+h,T-h]}(t)\\
\end{array}\right)_{t\in\tau_h}.	
\end{align*}
The continuous mapping theorem applied to (\ref{konvergenz_z_schlange_gemeinsam}) with map  $\varphi$ yields in  $(D[h,T-h],d_{SK})$ for $n\to\infty$ 
\begin{align*}
\varphi\left(\left(\widetilde Z_{1,t}^{(n)} \right)_{t\ge 0}, \left(\widetilde Z_{2,t}^{(n)} \right)_{t\ge 0}\right) \stackrel{d}{\longrightarrow}
\varphi\left(\left(\frac{\sqrt{\sigma_1^2/\mu_1^3}}{s_{t}^{(1)}} W_{1,t} \right)_{t\ge 0},\left(\frac{\sqrt{\sigma_2^2/\mu_2^3}}{s_{t}^{(1)}} W_{2,t} \right)_{t\ge 0} \right).
\end{align*}
Thus, it remains to be shown that 
\begin{align}
&\left(\Gamma_t^{(n)}\right)_{t\in\tau_h}= \varphi\left(\left(\widetilde Z_{1,t}^{(n)} \right)_{t\ge 0}, \left(\widetilde Z_{2,t}^{(n)} \right)_{t\ge 0}\right),\label{konvergenz_alternative_a}\\
&\left(L_t\right)_{t\in\tau_h}\sim \varphi\left(\left(\frac{\sqrt{\sigma_1^2/\mu_1^3}}{s_{t}^{(1)}} W_{1,t} \right)_{t\ge 0},\left(\frac{\sqrt{\sigma_2^2/\mu_2^3}}{s_{t}^{(1)}} W_{2,t} \right)_{t\ge 0} \right)\label{konvergenz_alternative_b},
\end{align}
where $\sim$ denotes equality in distribution. In order to show (\ref{konvergenz_alternative_a}) and (\ref{konvergenz_alternative_b}) we differentiate the four cases $t\in[h,c-h)$, $ t\in [c-h, c)$, $t \in [c, c+h)$ and $t\in [c+h,T-h]$.

%%%%
\noindent\underline{Derivation of (\ref{konvergenz_alternative_a}):}\\
Case $t<c-h:$
\begin{align*}
\left.\varphi\left(\left(\widetilde Z_{1,t}^{(n)} \right)_{t\ge 0}, \left(\widetilde Z_{2,t}^{(n)} \right)_{t\ge 0}\right)\right|_t
%& = \left(\frac{N_{1,n(t+h)}-\nicefrac{n(t+h)}{\mu_1}}{s_{h,t}^{(n)}}- \frac{N_{1,nt}-\nicefrac{nt}{\mu_1}}{s_{h,t}^{(n)}}\right) -  \left( \frac{N_{1,nt}-\nicefrac{nt}{\mu_1}}{s_{h,t}^{(n)}}-\frac{N_{1,n(t-h)}-\nicefrac{n(t-h)}{\mu_1}}{s_{h,t}^{(n)}}\right)\\
& = \frac{(N_{1,n(t+h)}-N_{1,nt}) -(N_{1,nt}-N_{1,n(t-h)})}{s_{t}^{(n)}}\\
& = \frac{[(N_{n(t+h)}^{(n)}-N_{nt}^{(n)}) -(N_{nt}^{(n)}-N_{n(t-h)}^{(n)})]-m_{t}^{(n)}}{s_{t}^{(n)}} = \Gamma_{t}^{(n)}.
\end{align*}
For $t\ge c+h$ we obtain analogous results by exchanging  subscripts. For $t\in [c-h,c)$ we obtain
\begin{align*}
&\left. \varphi\left(\left(\widetilde Z_{1,t}^{(n)} \right)_{t\ge 0}, \left(\widetilde Z_{2,t}^{(n)} \right)_{t\ge 0}\right)\right|_t\\
%& = \left(\frac{N_{2,n(t+h)}-\nicefrac{n(t+h)}{\mu_2}}{s_{h,t}^{(n)}}- \frac{N_{2,nc}-\nicefrac{nc}{\mu_2}}{s_{h,t}^{(n)}}\right) +  \left( \frac{N_{1,nc}-\nicefrac{nc}{\mu_1}}{s_{h,t}^{(n)}}-\frac{N_{1,nt}-\nicefrac{nt}{\mu_1}}{s_{h,t}^{(n)}}\right)\\
%& \qquad\qquad\qquad -   \left( \frac{N_{1,nt}-\nicefrac{nt}{\mu_1}}{s_{h,t}^{(n)}}-\frac{N_{1,n(t-h)}-\nicefrac{n(t-h)}{\mu_1}}{s_{h,t}^{(n)}}\right)\\
& = \frac{(N_{2,n(t+h)}-N_{2,nc}) + (N_{1,nc}-N_{1,nt}) -(N_{1,nt}-N_{1,n(t-h)}) - n \left(\frac{(t+h)-c}{\mu_2}-\frac{(t+h)-c}{\mu_1}\right)}{s_{t}^{(n)}}\\
& = \frac{[(N_{n(t+h)}^{(n)}-N_{nt}^{(n)}) -(N_{nt}^{(n)}-N_{n(t-h)}^{(n)})]-m_{t}^{(n)}}{s_{t}^{(n)}} = \Gamma_{t}^{(n)}.
\end{align*}
Analogously, we obtain $c\le t < c+h$, which proves (\ref{konvergenz_alternative_a}).

\noindent\underline{Derivation of  (\ref{konvergenz_alternative_b}):}\\
For $t< c-h$ we obtain
\begin{equation}\label{ko_al_b1}
\left.\varphi\left(\left(\frac{\sqrt{\sigma_1^2/\mu_1^3}}{s_{t}^{(1)}} W_{1,t} \right)_{t\ge 0},\left(\frac{\sqrt{\sigma_2^2/\mu_2^3}}{s_{t}^{(1)}} W_{2,t} \right)_{t\ge 0} \right)\right|_t
%& = \left( \frac{\sqrt{\nicefrac{\sigma_1^2}{\mu_1^3}}}{sd_{h,t}^{(1)}} W_{1,t+h} - \frac{\sqrt{\nicefrac{\sigma_1^2}{\mu_1^3}}}{sd_{h,t}^{(1)}} W_{1,t}\right) - \left(\frac{\sqrt{\nicefrac{\sigma_1^2}{\mu_1^3}}}{sd_{h,t}^{(1)}} W_{1,t} - \frac{\sqrt{\nicefrac{\sigma_1^2}{\mu_1^3}}}{sd_{h,t}^{(1)}} W_{1,t-h}\right)\nonumber\\
%& = \frac{\sqrt{\nicefrac{\sigma_1^2}{\mu_1^3}}}{\sqrt{2h\sigma_1^2/\mu_1^3}} \left[ \left(W_{1,t+h}-W_{1,t}\right)-\left(W_{1,t}-W_{1,t-h}\right)\right]\nonumber\\
= \frac{\left(W_{1,t+h}-W_{1,t}\right)-\left(W_{1,t}-W_{1,t-h}\right)}{\sqrt{2h}} = L_{t}.
\end{equation}
The same  holds for $t\ge c+h$ with the subscript exchanged. In the case $c-h\le t < c$ we obtain
\begin{align}\label{ko_al_b2}
& \left.\varphi\left(\left(\frac{\sqrt{\sigma_1^2/\mu_1^3}}{s_{t}^{(1)}} W_{1,t} \right)_{t\ge 0},\left(\frac{\sqrt{\sigma_2^2/\mu_2^3}}{s_{t}^{(1)}} W_{2,t} \right)_{t\ge 0} \right)\right|_t\nonumber\\
%& = \left( \frac{\sqrt{\nicefrac{\sigma_2^2}{\mu_2^3}}}{sd_{h,t}^{(1)}} W_{2,t+h} - \frac{\sqrt{\nicefrac{\sigma_2^2}{\mu_2^3}}}{sd_{h,t}^{(1)}} W_{2,c}\right)
%+  \left( \frac{\sqrt{\nicefrac{\sigma_1^2}{\mu_1^3}}}{sd_{h,t}^{(1)}} W_{1,c} - \frac{\sqrt{\nicefrac{\sigma_1^2}{\mu_1^3}}}{sd_{h,t}^{(1)}} W_{1,t}\right) \nonumber\\
%& \qquad\qquad\qquad - \left(\frac{\sqrt{\nicefrac{\sigma_1^2}{\mu_1^3}}}{sd_{h,t}^{(1)}} W_{1,t} - \frac{\sqrt{\nicefrac{\sigma_1^2}{\mu_1^3}}}{sd_{h,t}^{(1)}} W_{1,t-h}\right)\nonumber\\
& = \frac{\sqrt{\sigma_2^2/\mu_2^3} \left( W_{2,t+h} -  W_{2,c}\right) + \sqrt{\sigma_1^2/\mu_1^3} \left[ \left( W_{1,c}- W_{1,t}\right) - \left( W_{1,t}- W_{1,t-h}\right) \right] }{s_{t}^{(1)}}= L_{t}.
\end{align}
Analogously, we obtain for $c\le t < c+h$
\begin{align}\label{ko_al_b3}
&\left. \varphi\left(\left(\frac{\sqrt{\sigma_1^2/\mu_1^3}}{s_{t}^{(1)}} W_{1,t} \right)_{t\ge 0},\left(\frac{\sqrt{\sigma_2^2/\mu_2^3}}{s_{t}^{(1)}} W_{2,t} \right)_{t\ge 0} \right)\right|_t\nonumber\\
& = \frac{\sqrt{\sigma_2^2/\mu_2^3} \left[\left( W_{2,t+h} -  W_{2,t}\right) -\left( W_{2,t} -  W_{2,c}\right)\right] + \sqrt{\sigma_1^2/\mu_1^3}  \left( W_{1,c}- W_{1,t-h}\right) }{s_{t}^{(1)}}= L_{t}.
\end{align}
Now let $(W_t)_{t\ge0}$ be a standard Brownian motion, i.e., $(W_t)_{t\ge0} \sim (W_{1,t})_{t\ge0} \sim (W_{2,t})_{t\ge0}$. The process defined in (\ref{ko_al_b1}), (\ref{ko_al_b2}) and (\ref{ko_al_b3}) has continuous sample paths and is given as a function of increments of disjoint intervals of the processes $(W_{1,t})_{t\ge0}$ and $(W_{2,t})_{t\ge0}$. Therefore, we can omit the subscripts one and two in (\ref{ko_al_b1}), (\ref{ko_al_b2}) and (\ref{ko_al_b3}) and obtain a process that has continuous sample paths and the same distribution as the former one. By omitting the subscripts, we obtain the limit process $L$ as defined in equation (\ref{limit_process_cp}), which completes the proof of Proposition \ref{hauptaussage_alternative}. \hfill$\Box$\\
%Note that Proposition \ref{hauptaussage_alternative} refers to RPs, but that the proof applies directly to a larger class of renewal processes with a certain degree of variability in the variance defined in \citet{Messer2013}. 

\subsection{Proof of Lemma \ref{conv_delta}}\label{proof:conv_delta}

The Proof of Lemma \ref{conv_delta} works as follows: The uniform convergence $(s_t^{(n)}/\hat s_t^{(n)})_{t\in\tau_h}\to (\Delta)_{t\in\tau_h}$ a.s.~is equivalent to the uniform convergence $(\tilde s_t^{(n)}/\hat s_t^{(n)})_{t\in\tau_h}\to (1)_{t\in\tau_h}$ a.s.~as $n\to\infty$. The terms $\tilde s_t^{(n)}$ and $\hat s_t^{(n)}$ are functions of the estimators $\hat\mu_{le},\hat\mu_{ri},\hat\sigma_{le}^2$ and $\hat\sigma_{ri}^2$ as given in (\ref{schaetzer_s}). We show the uniform a.s.~convergence to their counterparts $\mu_{le},\mu_{ri},\sigma_{le}^2$ and $\sigma_{ri}^2$ defined in (\ref{mu_ri_cp}) and (\ref{sig_ri_cp}). More precisely, we show the uniform a.s.~convergence of $(\hat{\mu}_{le})_{t\in\tau_h}$ to $(\mu_{le})_{t\in\tau_h}$ and $(\hat{\mu}_{ri})_{t\in\tau_h}$ to $(\mu_{ri})_{t\in\tau_h}$ in Lemma \ref{lemm_conv_mu_cp}, and the uniform a.s.~convergence of $(\hat{\sigma}^2_{le})_{t\in\tau_h}$ to $(\sigma_{le}^2)_{t\in\tau_h}$ and $(\hat{\sigma}^2_{ri})_{t\in\tau_h}$ to $(\sigma_{ri}^2)_{t\in\tau_h}$ in Lemma \ref{lemm_conv_sig_cp}. 
Thus, the assertion of the Proposition holds true by the structure of the estimator $\hat s^2$ in (\ref{schaetzer_s}) and the function $\tilde s^2$ in (\ref{schaetzer_stilde}) and 
because convergence of sums and products of c\`{a}dl\`{a}g-valued functions in supremum norm is preserved when the limits are constant. 
 \hfill $\Box$\\

For completeness of the proof, we  show the consistency of the estimators  $\hat\mu_{le}$ and $\hat\mu_{ri}$ in Lemma \ref{lemm_conv_mu_cp} and the consistency of  $\hat\sigma_{le}^2$ and $\hat\sigma_{ri}^2$ in Lemma \ref{lemm_conv_sig_cp}. For that we first show a functional version of the SLLN in the following Lemma.

\begin{lemm}\label{lemm_conv_nt_cp} 
	For the counting process $N_t^{(n)}$ that corresponds to the process $\Phi^{(n)}$, it holds in $(D[h,T-h],d_{\|\cdot\|})$ as $n\to\infty$ almost surely
	\begin{align}\label{conv_nt1_cp}
	\left(\frac{N_{n(t+h)}^{(n)}-N_{nt}^{(n)}}{nh}\right)_{t\in\tau_h} &\longrightarrow \left(\frac{1}{\mu_{ri}(h,t)}\right)_{t\in\tau_h},\\
	\label{conv_nt2_cp}
	\left(\frac{N_{nt}^{(n)}-N_{n(t-h)}^{(n)}}{nh}\right)_{t\in\tau_h} &\longrightarrow \left(\frac{1}{\mu_{le}(h,t)}\right)_{t\in\tau_h}.
	\end{align}
\end{lemm}

\textbf{Proof:}	
Outline: We show the convergence of the right window half as stated in (\ref{conv_nt1_cp}). The statement for the left window half follows analogously.\\ 
First, we show that for all $t\ge 0$ and all $h>0$ it holds almost surely as $n\to\infty$
\begin{align}\label{conv_nt1_cp_marginal}
\frac{N_{n(t+h)}^{(n)}-N_{nt}^{(n)}}{nh} \longrightarrow \frac{1}{\mu_{ri}(h,t)}.
\end{align}
Then, by a discretization argument this result is extended to hold true in  $(D[h,T-h],d_{\|\,\cdot\,\|})$, as stated in (\ref{conv_nt1_cp}).

\noindent\underline{Derivation of (\ref{conv_nt1_cp_marginal}):}\\
In order to show the convergence in (\ref{conv_nt1_cp_marginal}), we distinguish between three cases. First assume $t\le c-h$. Here, for all $n=1,2,\ldots$, the corresponding window $(nt,n(t+h)]$ lies left of the change point $nc$. Thus, the counting process $N_t^{(n)}$ completely refers to the first RP $\Phi_1(\mu_1,\sigma_1^2)$, i.e., $N_{nt}^{(n)}=N_{1,nt}$, while $(N_{1,t})_{t\ge 0}$ denotes the counting process associated with $\Phi_1$. Then, it can be shown that it holds almost surely for $n\to\infty$
\begin{align}\label{conv_nt1_cp_marginal_a}
\frac{N_{n(t+h)}^{(n)}-N_{nt}^{(n)}}{nh} = \frac{N_{1,n(t+h)}-N_{1,nt}}{nh} \longrightarrow \frac{1}{\mu_1}=\frac{1}{\mu_{ri}(h,t)},
\end{align}
compare e.g., \citet{Messer2013}. An analogous statement holds for $t>c$.

For $t\in(c-h,c]$, the right window half refers partially to $\Phi_1$ and $\Phi_2$. The section $(nt,nc]$ refers to $\Phi_1$ and the section $(nc,n(t+h)]$ corresponds to $\Phi_2$. Thus, we decompose $N_{n(t+h)}^{(n)}-N_{nt}^{(n)} = (N_{2,n(t+h)}-N_{2,nc}) + (N_{1,nc}-N_{1,nt})$. We obtain almost surely for $n\to\infty$
\begin{align}\label{conv_nt1_cp_marginal_b}
\frac{N_{n(t+h)}^{(n)}-N_{nt}^{(n)}}{nh} 	& = \frac{(N_{2,n(t+h)}-N_{2,nc}) + (N_{1,nc}-N_{1,nt})}{nh} \nonumber\\
& = \frac{t+h-c}{h}\;\; \frac{N_{2,n(t+h)}-N_{2,nc}}{n(t+h-c)} + \frac{c-t}{h}\;\;\frac{N_{1,nc}-N_{1,nt}}{n(c-t)}\nonumber\\
& \longrightarrow \frac{t+h-c}{h}\;\; \frac{1}{\mu_2} +  \frac{c-t}{h}\;\; \frac{1}{\mu_1} = \frac{1}{\mu_{ri}(h,t)}.
\end{align}
In total, the convergences (\ref{conv_nt1_cp_marginal_a}) - (\ref{conv_nt1_cp_marginal_b}) yield (\ref{conv_nt1_cp_marginal}).\\

\noindent\underline{Derivation of (\ref{conv_nt1_cp}):}\\
In order to show that also convergence in $(D[h,T-h],d_{\|\cdot\|})$ holds, we even show the convergence in (\ref{conv_nt1_cp}) on $[0,T-h]$.
It is sufficient to show that almost surely
\begin{align} \label{conv_nt4_cp}
\lim_{n\to\infty} \sup_{t\in[0,T-h]}  \frac{N_{n(t+h)}^{(n)}-N_{nt}^{(n)}}{nh/\mu_{ri}(h,t)}   \le 1 \qquad \text{and} \qquad \lim_{n\to\infty} \inf_{t\in[0,T-h]} \frac{N_{n(t+h)}^{(n)}-N_{nt}^{(n)}}{nh/\mu_{ri}(h,t)} \ge 1.
\end{align}	

\begin{figure}[htb]
	\centering
	\includegraphics[scale=0.65]{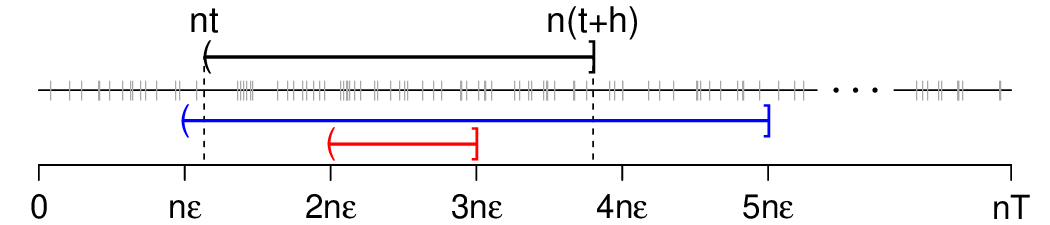}
	\vspace{0.3em}
	\caption{Schematic representation of the discretization of the time horizon $(0,nT]$ into equidistant sections of length $n\varepsilon$. All windows of length $nh$ (black) are contained in one of finitely many windows of length $n(\lceil h/\varepsilon \rceil + 1)\varepsilon$ (blue) and contain one of finitely many windows of length $n(\lfloor h/\varepsilon \rfloor - 1)\varepsilon$ (red).
	By letting $\varepsilon\downarrow0$, the size of the finitely many blue and red windows gets arbitrarily close to the true window size $nh$.}
	\label{alt}
\end{figure}

We show the left inequality of (\ref{conv_nt4_cp}). The right one follows analogously. 
We use a discretization argument. For $x\in\mathbb{R}$ let $|\!\lceil x \rceil\!|:= \lceil x \rceil + 1$.
For $\varepsilon>0$ with $T/\varepsilon\in\mathbb{N}$ we decompose the time interval $(0,nT]$ into equidistant sections of length $n\varepsilon$ (Figure \ref{alt}).
Then we observe a set $S_\varepsilon:=\{(kn\varepsilon,kn\varepsilon + n|\!\lceil h/\varepsilon \rceil\!| \varepsilon] :k=0,1,\ldots,T/\varepsilon - |\!\lceil h/\varepsilon \rceil\!|\}$  of finitely many windows of size $n|\!\lceil h/\varepsilon \rceil\!|\varepsilon$. The 
windows of $S_\varepsilon$ are slightly larger than $nh$ and for every $t\in(0,T-h]$ we find an element of $S_\varepsilon$ that overlaps the window $(nh,n(t+h)]$ (blue window in Figure \ref{alt}).
We bound 
\begin{align}\label{bound_sup_n_cp_3}
 \sup_{t\in[0,T-h]}  \frac{N_{n(t+h)}^{(n)}-N_{nt}^{(n)}}{nh/\mu_{ri}(h,t)} \le
& \max_{k=0,1,\ldots,T/\varepsilon - |\!\lceil h/\varepsilon \rceil\!|} \frac{N_{kn\varepsilon + n|\!\lceil h/\varepsilon \rceil\!| \varepsilon}^{(n)} - N_{kn\varepsilon}^{(n)}}{nh / \mu_{ri}(h,k\varepsilon)}.
\end{align} 
Now we make use of the fact that the convergence in (\ref{conv_nt1_cp}) holds true for a finite number of windows. By letting $n\to\infty$ the right hand side of (\ref{bound_sup_n_cp_3}) converges to $1+\delta_\varepsilon$, with
\begin{align*}
\delta_\varepsilon& \le   \frac{\max(\mu_1,\mu_2)}{\min(\mu_1,\mu_2)} \,\,\frac{|\!\lceil h/\varepsilon\rceil\!|\varepsilon-h}{h}. 
\end{align*}
The expression $\delta_\varepsilon>0$ accounts for the additional portion that results from the enlarged windows. Then, by letting $\varepsilon\downarrow 0$ the summand  $\delta_\varepsilon$ vanishes, which yields the first inequality in (\ref{conv_nt4_cp}). 
Analogously, for the lower bound of (\ref{conv_nt4_cp})  we find finitely many smaller windows of length $n(\lfloor h/\varepsilon \rfloor - 1)\varepsilon$, such that every window $(nh,n(t+h)])$ contains such a smaller window (red window in Figure \ref{alt}). 
Then, the limit of the infimum can be bounded from below by $1-\delta_\varepsilon'$ with $\delta_\varepsilon'>0$ and such that  $\delta_\varepsilon'\to 0$ as $\varepsilon\downarrow0$. Here, $\delta_\varepsilon'$ refers to the portion that is not covered by choosing the finitely many windows to be slightly smaller than the true window size $nh$.
\hfill$\Box$\\ 

Next, we show the uniform a.s.~convergences $(\hat\mu_{ri}(nh,nt))_{t\in\tau_h}\to(\mu_{ri}(h,t))_{t\in\tau_h}$ and\\ $(\hat\mu_{le}(nh,nt))_{t\in\tau_h}\to(\mu_{le}(h,t))_{t\in\tau_h}$ as $n\to\infty$. The estimators are given as
	\begin{align*}
		\hat\mu_{ri} =\hat\mu_{ri}(nh,nt)=\frac{1}{N_{n(t+h)}^{(n)}-N_{nt}^{(n)}-1}\sum_{i=N_{nt}^{(n)}+2}^{N_{n(t+h)}^{(n)}}\xi_i^{(n)}, \quad  \textrm{if}\quad N_{n(t+h)}^{(n)}-N_{nt}^{(n)}>1,
	\end{align*}
	and $\hat\mu_{ri}=0$ otherwise and $\hat\mu_{le}$ is given analogously.
	
	\begin{lemm}\label{lemm_conv_mu_cp}
		 For the estimators $\hat\mu_{ri}(nh,nt)$ and $\hat\mu_{le}(nh,nt)$ as given in (\ref{schaetzer_s}), it holds in $(D[h,T-h],d_{\|\cdot\|})$ as $n\to\infty$ almost surely
		\begin{align}
			(\hat\mu_{ri}(nh,nt))_{t\in\tau_h} \longrightarrow (\mu_{ri}(h,t))_{t\in\tau_h},\label{conv_mu1_cp}\\
			(\hat\mu_{le}(nh,nt))_{t\in\tau_h} \longrightarrow (\mu_{le}(h,t))_{t\in\tau_h}.
		\end{align}
	\end{lemm}
	
	\textbf{Proof:}
	 We show the convergence of the right window half as stated in (\ref{conv_mu1_cp}). The assertion for the left window half follows analogously. We proceed as in the proof Lemma \ref{lemm_conv_nt_cp}. First, we show that for all $t\ge 0$ and $h>0$, it holds almost surely as $n\to\infty$
	 \begin{align}\label{conv_mu1_cp_marginal}
	 \frac{1}{nh}\sum_{i=N_{nt}^{(n)}+2}^{N_{n(t+h)}^{(n)}}\xi_i^{(n)}\longrightarrow 1,
	 \end{align}
	 i.e., the sum of the life times in the window half  asymptotically equals the window length.
	 Then, this result is extended to $(D[h,T-h],d_{\|\,\cdot\,\|})$ to conclude the convergence in (\ref{conv_mu1_cp}). 
	 \noindent\underline{Derivation of (\ref{conv_mu1_cp_marginal}):}\\
	 Assertion (\ref{conv_mu1_cp}) has been shown in \citet{Messer2013}  to hold for the individual processes $\Phi_j(\mu_j,\sigma_j^2)$. Therefore, as we show (\ref{conv_mu1_cp}) for the right window half, convergence (\ref{conv_mu1_cp_marginal})
	 holds true for $t\in(0,c-h]$ and $t\ge c$. For $t\in(c-h,c]$, the right window half contains parts of $\Phi_1$ and of $\Phi_2$. We therefore decompose $(nt,n(t+h)] = (nt,nc] \cup (nc,n(t+h)]$. The section $(nt,nc]$ refers to $\Phi_1$ and the section $(nc,n(t+h)]$ corresponds to $\Phi_2$. The life time at the change point $c$ results from $\Phi_1$ and $\Phi_2$, and we therefore bound 
	 \begin{align}\label{bound_xi_cp}
	 \sum_{i=N_{1,nt}+2}^{N_{1,nc}}\xi_{1,i} + \sum_{i=N_{2,nc}+2}^{N_{2,n(t+h)}}\xi_{2,i}\le \sum_{i=N_{nt}^{(n)}+2}^{N_{n(t+h)}^{(n)}}\xi_i^{(n)} \le \sum_{i=N_{1,nt}+2}^{N_{1,nc}+1}\xi_{1,i} + \sum_{i=N_{2,nc}+1}^{N_{2,n(t+h)}}\xi_{2,i},
	 \end{align}
	which allows to use the properties of the individual processes. For the right hand side in (\ref{bound_xi_cp}) it holds as $n\to\infty$
	 \begin{align*}%\label{conv_mu1_cp_marginal_b}
	 \frac{1}{nh}\sum_{i=N_{nt}^{(n)}+2}^{N_{n(t+h)}^{(n)}}\xi_i^{(n)} 
	& \le \left(\frac{c-t}{h} \,\,\,\frac{1}{n(c-t)} \sum_{i=N_{1,nt}+2}^{N_{1,nc}+1}\xi_{1,i} \right)
	 +\left(\frac{t+h-c}{h} \,\,\,\frac{1}{n(t+h-c)} \sum_{i=N_{2,nc}+1}^{N_{2,n(t+h)}}\xi_{2,i}\right)\nonumber\\
	 &\longrightarrow \frac{(c-t)}{h} +\frac{t+h-c}{h}\,\, = 1.
	 \end{align*}
	 Analogously we obtain the lower bound, such that assertion (\ref{conv_mu1_cp_marginal}) holds true.
	 
	 \noindent\underline{Derivation of (\ref{conv_mu1_cp}):}\\
     In order to show convergence (\ref{conv_mu1_cp_marginal}) in  $(D[h,T-h],d_{\|\,\cdot\,\|})$ we even prove it on the interval  $[0,T-h]$. For that, we show that almost surely
	 \begin{align}\label{uniforme_konvergenz_cp_xi_2}
	 \lim_{n\to\infty} \sup_{t\in[0,T-h]}  \frac{1}{nh} \sum_{i=N_{nt}^{(n)}+2}^{N_{n(t+h)}^{(n)}}\xi_i^{(n)}  \,\le\, 1  \qquad \textrm{and} \qquad \lim_{n\to\infty} \inf_{t\in[0,T-h]} \frac{1}{nh} \sum_{i=N_{nt}^{(n)}+2}^{N_{n(t+h)}^{(n)}}\xi_i^{(n)} \ge 1.
	 \end{align}
	 We show here the left inequality of (\ref{uniforme_konvergenz_cp_xi_2}). We use the same discretization argument as in the proof of Lemma \ref{lemm_conv_nt_cp} and decompose the interval $(0,nT]$ into equidistant sections of length $n \varepsilon$ (Figure \ref{alt}).
	 Then we bound
	 \begin{align*}
	 \sup_{t\in [0,T-h]}  \frac{1}{nh} \sum_{i=N_{nt}^{(n)}+2}^{N_{n(t+h)}^{(n)}}\xi_i^{(n)} &\le \max_{k=0,1,\ldots,T/\varepsilon - |\!\lceil h/\varepsilon \rceil\!|}  \frac{1}{nh} \sum_{i=N_{kn\varepsilon}^{(n)}}^{N_{kn\varepsilon+n|\!\lceil h/\varepsilon\rceil\!|\varepsilon }^{(n)}}\xi_i^{(n)}\\
	 %& \le \max_{k=0,1,\ldots,T/\varepsilon - |\!\lceil h/\varepsilon \rceil\!|}  \frac{1}{nh} \sum_{i=N_{kn\varepsilon+nh}^{(n)}}^{N_{kn\varepsilon+n|\!\lceil h/\varepsilon\rceil\!|\varepsilon }^{(n)}}\xi_i^{(n)}+ \max_{k=0,1,\ldots,T/\varepsilon - |\!\lceil h/\varepsilon \rceil\!|}  \frac{1}{nh} \sum_{i=N_{kn\varepsilon}}^{N_{kn\varepsilon+nh }^{(n)}}\xi_i^{(n)}\\
	 & \le \frac{|\!\lceil h/\varepsilon\rceil\!| \varepsilon -h}{h}
	 + \max_{k=0,1,\ldots,T/\varepsilon - |\!\lceil h/\varepsilon \rceil\!|}  \frac{1}{nh} \sum_{i=N_{kn\varepsilon}^{(n)}}^{N_{kn\varepsilon+nh}^{(n)}}\xi_i^{(n)}.
	 \end{align*}
	 The first summand tends to zero as $\varepsilon\downarrow 0$ and is independent of $n$. Further, for every $\varepsilon>0$, the second summand converges to unity almost surely as $n \to\infty$, according to equation (\ref{conv_mu1_cp_marginal}).
	 Thus, the first inequality in (\ref{uniforme_konvergenz_cp_xi_2}) holds. The second inequality in (\ref{uniforme_konvergenz_cp_xi_2}) can be shown similarly.
	 Thus, convergence (\ref{conv_mu1_cp_marginal}) holds in $(D[h,T-h],d_{\|\,\cdot\,\|})$,  and together with Lemma \ref{lemm_conv_nt_cp} the assertion (\ref{conv_mu1_cp}) holds true by Slutsky's Theorem. \hfill$\Box$\\
	 
	To finish the proof of Lemma \ref{conv_delta}, we need to show the uniform a.s.~convergences  $(\hat\sigma_{ri}^2(nh,nt))_{t\in\tau_h}\to(\sigma_{ri}^2(h,t))_{t\in\tau_h}$ and $(\hat\sigma_{ri}^2(nh,nt))_{t\in\tau_h}\to(\sigma_{ri}^2(h,t))_{t\in\tau_h}$ as $n\to\infty$. The estimator $\hat\sigma_{ri}^2$ is given as 
	\begin{align*}
		\hat\sigma^2_{ri}(nh,nt) =\frac{1}{N_{n(t+h)}^{(n)}-N_{nt}^{(n)}-2}\sum_{i=N_{nt}^{(n)}+2}^{N_{n(t+h)}^{(n)}}\left(\xi_i^{(n)}-\hat\mu(nh,nt)\right)^2 , \quad  \textrm{if}\quad N_{n(t+h)}^{(n)}-N_{nt}^{(n)}>2,
	\end{align*}
	and $\hat\sigma^2_{ri}=0$ otherwise. Similarly $\hat\sigma^2_{le}$ is given. 
	\begin{lemm}\label{lemm_conv_sig_cp}
		For the estimators $\hat\sigma_{ri}^2(nh,nt)$ and $\hat\sigma_{le}^2(nh,nt)$ as given in (\ref{schaetzer_s}) it holds in $(D[h,T-h],d_{\|\cdot\|})$ as $n\to\infty$ almost surely
		\begin{align}
			(\hat\sigma_{ri}^2(nh,nt))_{t\in\tau_h} \longrightarrow (\sigma_{ri}^2(h,t))_{t\in\tau_h},\label{conv_sig1_cp}\\
			(\hat\sigma_{le}^2(nh,nt))_{t\in\tau_h} \longrightarrow (\sigma_{le}^2(h,t))_{t\in\tau_h}.
		\end{align}
	\end{lemm}
	
\textbf{Proof:}	
	Again we show the convergence of the right window half as given in  (\ref{conv_sig1_cp}), while the statement for the left window half follows analogously.\\ 
	First, we show that for all $t\ge 0$ and $h>0$, it holds almost surely as $n\to\infty$
	\begin{align}\label{conv_sig1_cp_marginal}
	\frac{\mu_{ri}(h,t)}{nh}\sum_{i=N_{nt}^{(n)}+2}^{N_{n(t+h)}^{(n)}}\left(\xi_i^{(n)}-\hat\mu_{ri}(nh,nt)\right)^2\longrightarrow \sigma_{ri}^2(h,t).
	\end{align}
Then, this result is extended to $(D[h,T-h],d_{\|\,\cdot\,\|})$  which yields (\ref{conv_sig1_cp}).
	
	\noindent\underline{Derivation of (\ref{conv_sig1_cp_marginal}):}\\
	In the following, let $\hat\mu_{j,ri}(nh,nt)$ denote the estimator that corresponds to $\Phi_j(\mu_j,\sigma_j^2)$. As before, $\hat\mu_{ri}(nh,nt)$ denotes the estimator that refers to the compound process $\Phi^{(n)}$.
	
	Note that (\ref{conv_sig1_cp_marginal}) was shown in \citet{Messer2013} to hold for the individual processes $\Phi_j(\mu_j,\sigma_j^2)$. Therefore, as we show (\ref{conv_sig1_cp}) for the right window, (\ref{conv_sig1_cp_marginal}) holds for $t\in(0,c-h]$ and $t\ge c$. For the remaining case $t\in(c-h,c]$, we recall that the right window half partially corresponds to $\Phi_1$ and $\Phi_2$. Again, we decompose $(nt,n(t+h)] = (nt,nc] \cup (nc,n(t+h)]$, where the sections $(nt,nc]$ and $(nc,n(t+h)]$ refer to $\Phi_1$ and $\Phi_2$, respectively. We decompose 
	\begin{align}\label{conv_sig1_cp_marginal_1}
	&\frac{\mu_{ri}(h,t)}{nh}\sum_{i=N_{nt}^{(n)}+2}^{N_{n(t+h)}^{(n)}}\left(\xi_i^{(n)}-\hat\mu_{ri}(nh,nt)\right)^2\nonumber\\
	& \qquad= \left(\frac{(c-t)\mu_{ri}(h,t)}{h\mu_1}\,\,\frac{\mu_1}{n(c-t)} \sum_{i=N_{1,nt}^{(n)}+2}^{N_{1,nc}^{(n)}}\left(\xi_{1,i}-\hat\mu_{ri}(nh,nt)\right)^2\right)\nonumber\\
	&\qquad\qquad +\left(\frac{(t+h-c)\mu_{ri}(h,t)}{h\mu_2}\,\,\frac{\mu_2}{n(t+h-c)} \sum_{i=N_{2,nc}^{(n)}+2}^{N_{2,n(t+h)}^{(n)}}\left(\xi_{2,i}-\hat\mu_{ri}(nh,nt)\right)^2\right) + o_{a.s.}(1)
	\end{align}
	The term $o_{a.s.}(1)$ accounts for the summand that corresponds to the single life time $\xi_{N_{nc}^{(n)}+1}^{(n)}$ that overlaps the change point and that is not respected in the first two terms of (\ref{conv_sig1_cp_marginal_1}). By Borel-Cantelli Lemma, the sequence $\{(\xi_{N_{nc}^{(n)}+1}^{(n)}-\hat\mu_{ri}(nh,nt))^2/nh\}_{n=1,2,\ldots}$ can be shown to vanish almost surely for $n\to\infty$ and is therefore abbreviated with $o_{a.s.}(1)$. For the first summand, we find almost surely as $n\to\infty$
	\begin{align}
	& \frac{\mu_1}{n(c-t)} \sum_{i=N_{1,nt}^{(n)}+2}^{N_{1,nc}^{(n)}}\left(\xi_{1,i}-\hat\mu_{ri}(nh,nt)\right)^2\nonumber\\
	& \qquad= \frac{\mu_1}{n(c-t)}\sum_{i=N_{1,nt}^{(n)}+2}^{N_{1,nc}^{(n)}}\left(\left[\xi_{1,i}-\hat\mu_{1,ri}(n(c-t),nt)\right] \, + \, [\hat\mu_{1,ri}(n(c-t),nt)-\hat\mu_{ri}(nh,nt)]\right)^2\nonumber\\
	& \qquad=\frac{\mu_1}{n(c-t)}\sum_{i=N_{1,nt}^{(n)}+2}^{N_{1,nc}^{(n)}} [\xi_{1,i}-\hat\mu_{1,ri}(n(c-t),nt)]^2\label{deco_sig_cp}\\
	& \qquad\qquad+ 2[\hat\mu_{1,ri}(n(c-t),nt)-\hat\mu_{ri}(nh,nt)]\left(\frac{\mu_1}{n(c-t)}\sum_{i=N_{1,nt}^{(n)}+2}^{N_{1,nc}^{(n)}} [\xi_{1,i}-\hat\mu_{1,ri}(n(c-t),nt)]\right)\nonumber\\
	& \qquad\qquad+ \frac{\mu_1}{n(c-t)}\left(N_{1,nc}^{(n)}-N_{1,nt}^{(n)}-1\right) [\hat\mu_{1,ri}(n(c-t),nt)-\hat\mu_{ri}(nh,nt)]^2\nonumber\\
	& \qquad\longrightarrow \sigma_1^2 + (\mu_1-\mu_{ri}(h,t))^2.\nonumber
	\end{align}
	The first summand in (\ref{deco_sig_cp}) shows the a.s.~convergence to $\sigma_1^2$ because it refers only to $\Phi_1$. The second summand in (\ref{deco_sig_cp}) vanishes a.s.~since the left term converges a.s.~according to Lemma \ref{lemm_conv_mu_cp} and the right term tends to zero a.s.~according to Lemmas \ref{lemm_conv_nt_cp} and  \ref{lemm_conv_mu_cp}.  The third summand in (\ref{deco_sig_cp}) tends to $(\mu_1-\mu_{ri}(h,t))^2$ a.s., since the term in the squared brackets converges to  $(\mu_1-\mu_{ri}(h,t))^2$ a.s.~according to Lemma \ref{lemm_conv_mu_cp}, while the scaled counting process converges to unity a.s.~due to Lemma \ref{lemm_conv_nt_cp}. 
	
	An analogous result can be obtained for the second summand of  (\ref{conv_sig1_cp_marginal_1}) which yields almost surely for $n\to\infty$	
	\begin{align}
	\frac{\mu_{ri}(h,t)}{nh}\sum_{i=N_{nt}^{(n)}+2}^{N_{n(t+h)}^{(n)}}\left(\xi_i^{(n)}-\hat\mu_{ri}(nh,nt)\right)^2
	&\longrightarrow \left(\frac{(c-t)\mu_{ri}(h,t)}{h\mu_1}\,\,[\sigma_1^2 + (\mu_1-\mu_{ri}(h,t))^2]\right)\nonumber\\%\label{conv_sig1_cp_marginal_2}\\
	&\quad\qquad +\left(\frac{(t+h-c)\mu_{ri}(h,t)}{h\mu_2}\,\, [\sigma_2^2 + (\mu_2-\mu_{ri}(h,t))^2]\right),\nonumber
	%&\qquad = \sigma_{ri}^2(h,t).\nonumber
	\end{align} 
	and elementary calculations yield equality to $\sigma_{ri}^2(h,t)$.
	The convergence in (\ref{conv_sig1_cp_marginal}) can be concluded using an analogous discretization argument as in the proofs of Lemmas \ref{lemm_conv_mu_cp} and \ref{lemm_conv_nt_cp}, such that the assertion (\ref{conv_sig1_cp}) can be concluded.\hfill$\Box$\\

\end{document}